\documentclass[11pt]{amsart}
\usepackage{latexsym,amssymb,amsthm,amsmath}

\def\R{\mathbb{R}}
\def\H{\mathbb{H}}
\def\Hyp{\mathcal{H}}
\def\dt{\partial_t}
\def\du{\partial_u}
\def\dv{\partial_v}

\def\nablao{\stackrel{{\rm o}}{\nabla}}

\def\proof{{\sc Proof. }}

\newtheorem{contor}{1}[section]

\newtheorem{theorem}[contor]{Theorem}

\newtheorem{proposition}[contor]{Proposition}
\newtheorem{remark}[contor]{Remark}

\newcommand{\gata}{\hfill\hskip -1cm \rule{.5em}{.5em}}

\title{Constant Angle Surfaces in $\H^2\times \R$}

\author[F. Dillen]{Franki Dillen}
\address[F. Dillen]{Katholieke Universiteit Leuven\\ Departement
Wiskunde\\ Celestijnenlaan 200 B\\ B-3001 Leuven\\ Belgium}
\email[F. Dillen]{franki.dillen@wis.kuleuven.be}

\author[M. Munteanu]{Marian Ioan Munteanu}
\address[M. Munteanu]{University 'Al.I.Cuza' of Ia\c si\\
Faculty of Mathematics\\
Bd. Carol I, no.11\\
700506 Ia\c si\\
Romania}
\email[M. Munteanu]{munteanu2001@hotmail.com}

\date{\today }

\thanks{The second author was partially supported by Grant CEEX ET
n.5883/2006 ANCS (Romania) and Grant n.9/2006 in the framework of
the Bilateral Program between Italy and Romania.}

\begin{document}



\begin{abstract}
In last years, the study of the geometry of surfaces in product
spaces $S^2\times\R$ and $\H^2\times\R$ is developed by a large
number of mathematicians. In \cite{kn:DFVV07} the authors study
constant angle surfaces in $S^2\times \R$, namely those surfaces for
which the unit normal makes a constant angle with the tangent
direction to $\R$. In this paper we classify constant angle surfaces
in $\H^2\times\R$, where $\H^2$ is the hyperbolic
plane.
\end{abstract}



\maketitle

\section{Preliminaries}

Let $\widetilde M=\H^2\times\R$ be the Riemannian product of $\left(\H^2(-1),g_H\right)$ and $\R$ with the
standard Euclidean metric. Denote by $\widetilde g$ the product metric and by $\widetilde\nabla$
the Levi Civita connection of $\widetilde g$. Denote by $t$ the (global) coordinate on $\R$ and
hence $\dt=\frac\partial{\partial t}$ is the unit vector field in the tangent bundle
$T(\H^2\times\R)$ that is tangent to the $\R$-direction. Hence, the product metric can be
written as
$$
\widetilde g=g_H+dt^2.
$$

The Riemann-Christoffel curvature tensor $\widetilde R$ of $\H^2\times\R$ is given by
$$
\widetilde
R(X,Y,Z,W)=-g_H(X_H,W_H)g_H(Y_H,Z_H)+g_H(X_H,Z_H)g_H(Y_H,W_H)
$$
for any $X,Y,Z,W$ tangent to $\H^2\times\R$. If $X$ is a tangent
vector to $\H^2\times\R$ we put $X_H$ its projection to the tangent
space of $\H^2$.

Let $M$  be a surface in $\widetilde M=\H^2\times\R$. If $\xi$ is a
unit normal to $M$, then the shape operator is denoted by $A$. We
have the formulas of Gauss and Weingarten

{\bf (G)}\qquad\qquad $\widetilde\nabla_XY=\nabla_XY+h(X,Y)$

{\bf (W)}\qquad\qquad $\widetilde\nabla_X\xi=-AX$,

for every $X$ and $Y$ tangent to $M$. Here $\nabla$ is the Levi
Civita connection on $M$ and $h$ is a symmetric $(1,2)$-tensor field
taking values in the normal bundle and called the second fundamental
form of $M$. We have $\widetilde g(h(X,Y),\xi)=g(X,AY)$ for all
$X,Y$ tangent to $M$, where $g$ is the restriction of $\widetilde g$
to $M$.

Since $\dt$ is of unit length, we decompose $\dt$ as
\begin{equation}
\label{eq:dec_dt}
\dt=T+\cos\theta\ \xi
\end{equation}
where $T$ is the projection of $\dt$ on the tangent space of $M$ and
$\theta$ is the angle function,  defined by
\begin{equation}
\label{eq:theta}
\cos\theta=\widetilde g(\dt,\xi).
\end{equation}

If $X,Y$ are tangent to $M$, then we have the following relation
$$
g_H(X_H,Y_H)=g(X,Y)-g(X,T)g(Y,T).
$$
Thus, if $R$ is the Riemannian curvature on $M$ the equation of Gauss can be written as

\begin{tabbing}
{\bf (EG)} \quad $R(X,Y,Z,W)=$\= $\ g(AX,W)g(AY,Z)-g(AX,Z)g(AY,W)-
$\\[1mm]
\>
$
-g(X,W)g(Y,Z)+g(X,Z)g(Y,W)+
$\\[1mm]
\>
$
+g(X,W)g(Y,T)g(Z,T)+g(Y,Z)g(X,T)g(W,T)-
$\\[1mm]
\>
$-g(X,Z)g(Y,T)g(W,T)-g(Y,W)g(X,T)g(Z,T)$
\end{tabbing}
for every $X,Y,Z,W\in T(M)$.

Using the expression of the curvature $\widetilde R$ of $\H^2\times\R$, after a straightforward computation
we write the equation of Codazzi

{\bf (EC)} \qquad\quad $\nabla_XAY-\nabla_YAX-A[X,Y]=\cos\theta\big(g(X,T)Y-g(Y,T)X\big)$

for all $X,Y\in T(M)$.

Now we give the following
\begin{proposition}
Let $X$ be a tangent vector to $M$. We have
\begin{equation}\label{shape}
\left\{\begin{array}{l}
\nabla_XT=\cos\theta AX\\[2mm]
X(\cos\theta)=-g(AX,T).
\end{array}\right.
\end{equation}
\end{proposition}
\proof
For any $X$ tangent to $M$ we can write
$$
X=X_H+g(X,T)\ \dt\ .
$$
We have
$$
\widetilde\nabla_X\dt=\widetilde\nabla_{X_H}\dt+g(X,T)\widetilde\nabla_{\dt}\dt=0.
$$
On the other hand,
$$
\widetilde\nabla_X\dt=\widetilde\nabla_XT+\widetilde\nabla_X(\cos\theta)\xi=
    \nabla_XT+h(X,T)+X(\cos\theta)\xi-(\cos\theta)AX.
$$
Identifying the tangent and the normal part respectively, one gets
$$
\nabla_XT=\cos\theta AX\quad{\rm and}\quad X(\cos\theta)\xi=-h(X,T).
$$
Hence the conclusion.
\gata

\vspace{2mm}

If $\theta\in [0,\pi)$ is a constant angle, from the previous
proposition the following relation holds: $g(AX,T)=0$ for every $X$ tangent to
$M$ (at $p$), which is equivalent to
\begin{equation}
g(AT,X)=0\ ,\quad \forall X\in T_p(M).
\end{equation}
This means that, if $T\ne 0$, $T$ is a principal direction with
principal curvature $0$.

\begin{remark}
\rm If $T=0$ on $M$, then $\dt$ is always normal so, $M\subseteq
\H^2\times \{t_0\}$, for $t_0\in\R$.
\end{remark}

If $T\neq 0$, we consider
\begin{equation}
e_1=\frac 1 {||T||}\ T,
\end{equation}
where $||T||=\sin\theta$. \\[2mm]
Let $e_2$ be a unit vector tangent to $M$ and perpendicular to
$e_1$. Then the shape operator $A$ takes the following form
$$
S=\left(\begin{array}{cc}
 0 & 0 \\ 0 & \lambda
\end{array}\right)
$$
for a certain function $\lambda$ on $M$. Hence we have
\begin{equation}
\label{eq:h}
h(e_1,e_1)=0, \ h(e_1,e_2)=0,\ h(e_2,e_2)=\lambda\ \xi.
\end{equation}
\begin{proposition}
If $M$ is a constant angle surface in $\H^2\times\R$ with constant
angle $\theta\ne 0$, then $M$ has constant Gaussian curvature
$K=-\cos^2\theta$ and the projection $T$ of $\frac\partial{\partial
t}$ is a principal direction with principal curvature $0$.
\end{proposition}
\proof
We have to prove only the first part of this statement. To do this, we decompose
$e_1,e_2\in T(M)$ as
\begin{equation}
e_1=E_1+\sin\theta \dt\ ,\quad e_2=E_2
\end{equation}
with $E_1,E_2\in\chi(\H^2)$. We immediately have
$$
g_H(E_1,E_1)=\cos^2\theta, \ g_H(E_1, E_2)=0,\ g_H(E_2, E_2)=1.
$$
Combining (\ref{eq:h}) with Gauss' equation we find for Gaussian curvature of $M$
\begin{equation}
K=-\cos^2\theta.
\end{equation}
\gata

We conclude this section with the following
\begin{proposition}
The Levi Civita connection of $g$ on $M$ is given by
\begin{equation}
\label{eq:LCg}
\begin{array}{ll}
\nabla_{e_1}e_1=0, \ \nabla_{e_2}e_1=\lambda\cot\theta\ e_2, \\
\nabla_{e_1}e_2=0, \ \nabla_{e_2}e_2=-\lambda\cot\theta\ e_1.
\end{array}
\end{equation}
\end{proposition}
\proof Direct computation from (\ref{shape}). \gata

\section{Characterization of constant angle surfaces}

In this section we want to classify the constant angle surfaces $M$
in $\H^2\times\R$. There exist two trivial cases, namely $\theta=0$
and $\theta=\frac\pi2$. As we have already seen, in the first case
one has that $\frac\partial{\partial t}$ is always normal and hence
$M$ is an open part of $\H^2\times\{t_0\}$, while in the second case
$\frac\partial{\partial t}$ is always tangent. This corresponds to
the Riemannian product of a curve in $\H^2$ and $\R$.
\label{pg:trivial}
\\[2mm]
We can take coordinates $(u,v)$ on $M$ such that the metric $g$ on $M$ has the form
\begin{equation}
\label{eq:metric_g}
g=du^2+\beta^2(u,v)\ dv^2
\end{equation}
with
$\du:=\frac\partial{\partial u}=e_1$ and
$\dv:=\frac\partial{\partial v}=\beta\ e_2$,
where $\beta$ is a smooth function on $M$.
This can be done since $[e_1,e_2]$ is collinear with $e_2$.
We have
$$
0=[\du,\dv]=[\du,\beta
e_2]=\beta_ue_2+\beta[e_1,e_2]=\left(\beta_u-\beta\lambda\cot\theta\right)
e_2
$$
and hence $\beta$ satisfies the following PDE
\begin{equation}
\label{eq:PDEbeta}
\beta_u=\beta\lambda\cot\theta.
\end{equation}
We can now write the Levi Civita connection of $g$ on $M$ in terms
of coordinates $u$ and $v$, namely
\begin{equation}
\label{eq:LCg_uv}
\begin{array}{c}
\nabla_{\du}\du=0, \ \nabla_{\du}\dv=\nabla_{\dv}\du=\lambda\cot\theta\ \dv,\\[2mm]
\nabla_{\dv}\dv=-\beta\beta_u\du+\frac{\beta_v}\beta\ \dv.
\end{array}
\end{equation}
\begin{proposition}
The two functions $\lambda$ and $\beta$ are given by
\begin{equation}
\label{eq:lambda}
\lambda(u,v)=\sin\theta \tanh\left(u\cos\theta+C(v)\right)
\end{equation}
\begin{equation}
\label{eq:beta}
\beta(u,v)=D(v)\cosh\left(u\cos\theta+C(v)\right),
\end{equation}
or
\begin{equation}
\label{eq:lambda_const}
\lambda(u,v)=\pm\sin\theta
\end{equation}
\begin{equation}
\label{eq:beta_const}
\beta(u,v)=D(v) e^{\pm u\cos\theta}
\end{equation}
where $C$ and $D$ are smooth functions depending on $v$, $D(v)\neq0$ for any $v$.
\end{proposition}
\proof
From the equation of Codazzi (EC), if we put $X=e_1$ and $Y=e_2$ one obtains
that $\lambda$ must satisfy the following PDE
\begin{equation}
\lambda_u=\sin\theta\cos\theta-\lambda^2\cot\theta.
\end{equation}
By integration, one gets (\ref{eq:lambda}) or (\ref{eq:lambda_const}).
Now, solving (\ref{eq:PDEbeta}) we obtain $\beta$.

\gata

\vspace{2mm}

There are many models for the hyperbolic plane. In the following we will deal
with the Minkowski model or the hyperboloid model $\Hyp$ of $\H^2$, namely
the upper sheet $(z>0)$ of the hyperboloid
$$
   \{(x,y,z)\in {\R^3_1}\ :\ x^2+y^2-z^2=-1\}.
$$
We denoted by $\R^3_1$ the Minkowski 3-space with
Lorentzian metric tensor
$$<\cdot,\cdot>=dx^2+dy^2-dz^2.$$
The unit normal to $\Hyp$ in a point $p\in {\R^3_1}$ is $N=\pm p$.
We will take the external normal $N=p$ and we have
$<N,N>=-1$.

\vspace{2mm}

We recall the notion of the Lorentzian cross-product (see e.g.
\cite{kn:Rat06}):
$$
   \boxtimes:\R^3_1\times\R^3_1 \longrightarrow \R^3_1,
({}^t[a_1, a_2,a_3] , {}^t[b_1,b_2,b_3] )\mapsto
   {}^t[a_2b_3-a_3b_2, a_3b_1-a_1b_3,a_2b_1-a_1b_2].
$$

Let $M$ be a 2-dimensional surface in
$\Hyp\times\R\subset\R^3_1\times\R$. On the ambient space we
consider the product metric: $g_{\rm o}=dx^2+dy^2-dz^2+dt^2$. Denote
by $\nablao$ the Levi Civita connection on $\R^3_1\times\R$ and let
$D^\bot$ be the normal connection of $M$ in $\R^3_1\times\R$. If
$\widetilde\xi$ is the unit normal to $\widetilde M$, then
$\widetilde \xi(p_1,p_2,p_3,p_4)=(p_1,p_2,p_3,0)$.  The shape
operator w.r.t. $\widetilde\xi$ is denoted by $\widetilde A$.

\begin{theorem}
A surface $M$ in $\Hyp\times\R$ is a constant angle surface if and
only if the position vector $F$ is, up to isometries of
$\Hyp\times\R$, locally given by
\begin{equation}
\label{eq:immersionF_hyp}
F(u,v)=\big(\cosh(u\cos\theta)f(v)+\sinh(u\cos\theta)f(v)\boxtimes f'(v), u\sin\theta\big),
\end{equation}
where $f$ is a unit speed curve on $\Hyp$.
\end{theorem}
\begin{remark}\rm
This result is similar to that given in {\bf{Theorem 2}} in \cite{kn:DFVV07}.
\end{remark}

{\sc Proof of the Theorem}.

First we have to prove that the given immersion (\ref{eq:immersionF_hyp})
is a constant angle surface in $\Hyp\times\R$. To do this we compute the
tangent vectors (in an arbitrary point on $M$)
$$
    \begin{array}{rcl}
       F_u(u,v) & = & \big(\cos\theta\ [\sinh(u\cos\theta)f(v)+\cosh(u\cos\theta)f(v)\boxtimes f'(v)], \sin\theta\big)\\[2mm]
       F_v(u,v) & = & \big(\cosh(u\cos\theta)f'(v)+\sinh(u\cos\theta)f(v)\boxtimes f''(v),0\big)=\\[1mm]
                & = & \big(\ [\cosh(u\cos\theta)-\kappa(v)\sinh(u\cos\theta)]\ f'(v),0\big),
    \end{array}
$$
where $\kappa$ is the geodesic curvature of the curve $f$. This follows from the identity
$f''(v)=f(v)+\kappa(v) f(v)\boxtimes f'(v)$
which implies
$f(v)\boxtimes f''(v)=-\kappa(v)f'(v)$, since $f$ is a unit speed curve on $\Hyp$.

We will calculate now both $\xi$ and $\widetilde\xi$. The second normal vector is nothing
but the position vector where we take the last component to be 0, namely we have
$$
  \widetilde\xi (u,v)=\big( \cosh(u\cos\theta)f(v)+\sinh(u\cos\theta)f(v)\boxtimes f'(v),0\big).
$$
Looking for the expression of the unitary normal $\xi$ as linear combination of $f$, $f'$
and $f\boxtimes f'$ we find after some easy computations that
$$
   \xi(u,v)=\big(-\sin\theta\ [\sinh(u\cos\theta)f(v)+\cosh(u\cos\theta)f(v)\boxtimes f'(v)], \cos\theta\big).
$$
It follows $<\xi,\dt>=\cos\theta$ (which is a constant).

\vspace{2mm}

Conversely, consider a surface $M$ in $\Hyp\times\R$ with constant angle function $\theta$.
If $M$ is one of the trivial cases (see page \pageref{pg:trivial}), then it can be
parameterized by (\ref{eq:immersionF_hyp}). Suppose from now on that $\theta\notin\{0,\frac \pi2\}$.

We have

$$
(F_4)_u = \widetilde g(F_u, \partial_t) = \widetilde
g(F_u,T+\cos\theta \xi)=g(\du,T)=\sin\theta
$$
and
$$
(F_4)_v = \widetilde g(F_v, \partial_t) =g(\dv,T)=0.
$$
These relations and the initial condition $F_4(0,0)=0$ yield
\begin{equation}
\label{eq:F4}
F_4=u\sin\theta.
\end{equation}


If $X=(X_1,X_2,X_3,X_4)$ is tangent to $M$, then
$\nablao_X\widetilde\xi=(X_1,X_2,X_3,0)$. It follows
\begin{itemize}
\item[$\bullet$]  $ D^\bot_X\widetilde\xi=<(X_1,X_2,X_3,0),\xi>\xi=-\cos\theta<X,T>\xi$
\item[$\bullet$]  $D^\bot_X\xi=\cos\theta<X,T>\widetilde\xi$.
\end{itemize}

Then  $\widetilde A$ can be expressed in the basis $\{\du,\dv\}$ by
$$
\widetilde A=\left(\begin{array}{rr}
-\cos^2\theta & 0\\0 & -1
\end{array}\right).
$$

From (\ref{eq:dec_dt}) one has $\xi_j=-\tan \theta (F_j)_u$ for all
$j=1,2,3$.
 \vspace{2mm}

From the previous relations, the formula of Gauss can be written as:
\begin{equation}
\label{eq:35}
       (F_j)_{uu}=\cos^2\theta F_j
\end{equation}
\begin{equation}
\label{eq:36}
       (F_j)_{uv}=\lambda\cot\theta (F_j)_{v}
\end{equation}
\begin{equation}
\label{eq:37}
      (F_j)_{vv}=-\beta\beta_u (F_j)_{u}+\frac{\beta_v}{\beta}\ (F_j)_{v}-\lambda\beta^2\tan\theta (F_j)_{u}+
    \beta^2 F_j.
\end{equation}

Case 1: $\lambda$ satisfies (\ref{eq:lambda}).   Integrating
(\ref{eq:36}) one gets
$$
          (F_j)_{v}=H_j(v)\cosh(u\cos\theta+C(v)).
$$
Hence
$$
    F_j=\displaystyle\int\limits_0^v\cosh(u\cos\theta+C(\tau))H_j(\tau)d\tau+I_j(u).
$$
Substituting in (\ref{eq:35}) we obtain
$$
(*)\quad   I_j=K_j\cosh(u\cos\theta)+L_j\sinh(u\cos\theta)
$$
with $K_j$ and $L_j$ real constants.

We define the following functions

\qquad\qquad
$f_j=K_j+\displaystyle\int\limits_0^v\cosh C(\tau)H_j(\tau)d\tau$, ($j=1,2,3$)

\qquad\qquad
$g_j=L_j+\displaystyle\int\limits_0^v\sinh C(\tau)H_j(\tau)d\tau$, ($j=1,2,3$).

\vspace{2mm}

Case 2:  $\lambda$ satsifies (\ref{eq:lambda_const}). One gets
$$
F_j=e^{\pm u\cos\theta}\int\limits_0^vH_j(\tau)d\tau+I_j(u)
$$
with $I_j$ having the same form as in $(*)$.

In this case we put

\qquad\qquad
$f_j=K_j+\int\limits_0^vH_j(\tau)d\tau$, ($j=1,2,3$)

\qquad\qquad
$g_j=L_j\pm\int\limits_0^vH_j(\tau)d\tau$, ($j=1,2,3$).

\vspace{2mm}

Let $f=(f_1,f_2,f_3)$ and $g=(g_1, g_2, g_3)$.

\vspace{2mm}

To summarize, in both cases $F$ is of the following form:
\begin{equation}
        F=\big(\cosh(u\cos\theta) f+\sinh(u\cos\theta) g, u\sin\theta\big).
\end{equation}

\vspace{2mm}

Let $\epsilon_1=\epsilon_2=1$ and $\epsilon_3=-1$. Then we have
$$
 1. \sum \epsilon_j F_j^2=-1, \quad
 2. \sum \epsilon_j F_{j,u}^2=\cos^2\theta, \quad
 3. \sum \epsilon_j F_{j,u}F_{j,v}=0,\quad
 4. \sum \epsilon_j F_{j,v}^2=\beta^2.
$$

From 1. and 2. one obtains

\qquad\qquad\ $\sum \epsilon_j f^2_j - \sum \epsilon_j g^2_j=-2$.
\hfill{(i)}

Now, relations 1. and 3. can be written as

\qquad\qquad
$\begin{array}{l}
\sum \epsilon_j f^2_j \cosh^2(u\cos\theta)+\sum \epsilon_j g^2_j \sinh^2(u\cos\theta)+\\
        \qquad\qquad + 2\sum \epsilon_j f_jg_j \sinh(u\cos\theta)\cosh(u\cos\theta)=-1.
\end{array}
$ \hfill{(ii)}

and

\qquad\qquad
$\begin{array}{l}
\left(\sum \epsilon_j f_jf'_j+\sum \epsilon_j g_jg'_j\right)\sinh(u\cos\theta)\cosh(u\cos\theta)+\\
   \qquad\qquad +\sum \epsilon_j f'_jg_j\cosh^2(u\cos\theta)+\sum \epsilon_j f_jg'_j\sinh^2(u\cos\theta)=0.
   \end{array}
$ \hfill{(iii)}

By derivation in 1. one has

\qquad\qquad
$\begin{array}{l}
\sum \epsilon_j f_jf'_j \cosh^2(u\cos\theta)+\sum \epsilon_j g_jg'_j\sinh^2(u\cos\theta)+\\
    \qquad\qquad +\left(\sum \epsilon_j f_jg'_j+\sum \epsilon_j f'_jg_j\right) \sinh(u\cos\theta)\cosh(u\cos\theta)=0
 \end{array}
$    \hfill{(iv)}

\qquad\qquad
$\begin{array}{l}
\left(\sum \epsilon_j f_j^2+\sum \epsilon_j g_j^2\right)\sinh(u\cos\theta)\cosh(u\cos\theta)+\\
  \qquad\qquad  +\sum \epsilon_j f_jg_j (\cosh^2(u\cos\theta)+\sinh^2(u\cos\theta))=0.
  \end{array}
$       \hfill{(v)}

Finally (i), (ii) and (v) yield

\qquad\qquad\ $\sum \epsilon_j f_j^2=-1$,\quad  $\sum \epsilon_j
g_j^2=1$,\quad $\sum \epsilon_j f_jg_j=0$.

\vspace{2mm}

Moreover, it follows
\begin{equation*}
          \sum \epsilon_j f_jf'_j=0,\
        \sum \epsilon_j g_jg'_j=0,\
        \sum \epsilon_j f_jg'_j+\sum \epsilon_j f'_jg_j=0.
\end{equation*}

\vspace{2mm}

From (iii) we get
\begin{equation}
\sum \epsilon_j f'_jg_j=0\ {\rm and}\ \sum \epsilon_j f_jg'_j=0.
\end{equation}
Hence, the relation (iv) is identically satisfied.

\vspace{2mm}

Let's write these last equations in another way:

\begin{equation}
\begin{array}{l}
<f,f>=-1\\
<g,g>=1\\
<f,g>=0\\
<f',g>=0\\
<f,g'>=0
\end{array}
\quad {\rm and}
\quad
\begin{array}{l}
<f,f'>=0\\
<g,g'>=0
\end{array}
\end{equation}

\vspace{2mm}

We still have to develop the relation 4. This yields
$$
<H(v),H(v)>=<f',f'>-<g',g'>=D^2(v),\quad H=(H_1,H_2,H_3).
$$
Remark that $f$ can be thought as a curve on $\Hyp$ (while $g$ not).

\vspace{2mm}

Since $<f',f'>\geq 0$ (it can be easily proved), one can change the $v$-coordinate
such that $f$ becomes a unit speed curve in $\Hyp$; this corresponds to
$D(v)^2\cosh^2 C(v)=1$ or $D(v)^2=1$ (this depends of the value of $\lambda$).

\vspace{2mm}

We have $g\perp f$ and $g\perp f'$. This means that $g$ is collinear
to $f\boxtimes f'$. We have $<g,g>=1$ and $<f\boxtimes f',f\boxtimes
f'>=1$ and hence $g=\pm f\boxtimes f'$. We can assume that
$g=f\boxtimes f'$. Then $F$ is given by (\ref{eq:immersionF_hyp}) as
we wanted to prove. \gata

\begin{remark}\rm
Looking for all minimal constant angle surfaces in $\Hyp\times\R$,
these must be totally geodesic in $\Hyp\times\R$. Hence we obtain the following surfaces:

\qquad (1) $\Hyp\times\{t_0\}$, $t_0\in\R$

\qquad (2) $f\times\R$ with $f$ a geodesic line in $\Hyp$.
\end{remark}

\

\end{document}